\documentclass[a4paper,11pt]{amsart}

	\usepackage[utf8]{inputenc}
	\usepackage[T1]{fontenc}
	\usepackage[english]{babel}
    \usepackage{float}
	\usepackage[protrusion=true,expansion=true,kerning, spacing]{microtype}
	\microtypecontext{spacing=english}
	\usepackage{CJKutf8} % for Japanese and Kanji
	\usepackage{cmap} % PDF besser durchsuch- und kopierbar machen	

	\usepackage{xspace}
	
	\usepackage[textwidth=16cm,hcentering,heightrounded,foot=1cm]{geometry}
	\usepackage{caption}% for caption of figures outside float.
 \usepackage{titlesec}
	\titleformat{name=\section}[hang]{\scshape \large \centering}{\thesection.~}{0em}{}[]

	\titleformat{\subsection}[hang]{\bfseries}{\thesubsection.~}{0em}{}[]
	\titleformat{\subsubsection}[runin]{\bfseries}{\thesubsubsection.}{.2em}{}[]
 \usepackage{titletoc}
\titlecontents{section}[2pc]{\addvspace{1pc}}
{\contentslabel[\textsc{\thecontentslabel}]{1pc}}
{}{\dotfill\contentspage}
[\addvspace{2pt}]

\titlecontents{subsection}[4pc]{}
{\contentslabel[\subsectionname\thecontentslabel]{2pc}}
{}{\hfill\contentspage}
[]

\usepackage{csquotes}
\usepackage[% BIBLATEX MANUAL p.47
bibencoding=utf8,%
giveninits=true,
bibstyle=ieee-alphabetic,%ieee-alphabetic,% only option for dashed authors
citestyle=alphabetic,%
sorting=nyt,%
language=french,%
isbn=false,%
doi=false,%
eprint=true,%pour citation arxiv
backref=true, % back ref - at least for control during writing
hyperref=true,
url=true,
minbibnames=4,
maxbibnames=4,
dashed=true,
backend=biber,
]{biblatex}
	\usepackage{amsthm,amssymb,amsmath}
	\usepackage{amssymb,amsfonts}
	\usepackage[fixamsmath]{mathtools} %fixe un bug ams
	\usepackage{accents} % accents on math
	\usepackage{tikz-cd}
		\numberwithin{equation}{subsection}

	\usepackage{graphics,graphicx}
	\usepackage{color}
	\usepackage{tikz}
	\usepackage{caption} % for subfigure
	\usepackage{subcaption} % for subfigure
	\usepackage{wrapfig}

\theoremstyle{definition}
	\newtheorem{defin}{Definition}[section]

\theoremstyle{definition}
	\newtheorem{rem}[defin]{Remark}
	
\theoremstyle{plain}
	\newtheorem{theo}[defin]{Theorem}
	\newtheorem{prop}[defin]{Proposition}
	\newtheorem{cor}[defin]{Corollary}
	
	\newtheorem{quest}[defin]{Question}
	
\newcommand{\C}{\mathbf{C}}

\newcommand{\R}{\mathbf{R}}
\newcommand{\Q}{\mathbf{Q}}
\newcommand{\Z}{\mathbf{Z}}

\newcommand{\F}{\mathbf{F}}

\newcommand{\G}{\mathbf{G}}

\newcommand{\HH}{\mathbf{H}}

\newcommand{\OO} {{\mathcal O}} % \O already defined

\newcommand{\mfS} {\mathfrak{S}}

\newcommand{\Spec} { \operatorname{Spec} }

\newcommand{\End} {\operatorname{End} }

\newcommand{\Card} {\operatorname{Card} }
\newcommand{\Aut} {\operatorname{Aut} }
\newcommand{\Gal} {\operatorname{Gal} }

\newcommand{\GL}{\operatorname{GL}}

\newcommand{\Sp}{\operatorname{Sp}}

\usepackage[pdftex,
hidelinks,
hyperindex=true, bookmarks=true, bookmarksopen=true, bookmarksnumbered=true,
pdfauthor={Séverin Philip},
pdftitle={The local obstruction to semi-stable reduction for abelian varieties},
pdfsubject={},
pdfstartview=FitH
]{hyperref}

\begin{document}
\author{Séverin Philip}
\address{Department of Mathematics, Stockholms universitet
SE-106 91 Stockholm, Sweden}
	\email{severin.philip@math.su.se}
\thanks{This research was supported by the Verg foundation. \\
The author thanks J.~Bergström, B.~Collas and G.~Rémond for helpful comments on earlier versions of this text.}
\title[Local obstruction to semi-stable reduction]{The local obstruction to semi-stable reduction for abelian varieties \\[0.5em]
\textnormal{ Version of \today}}

\begin{abstract}
Grothendieck defined a group that represents the local obstruction for an abelian variety to have semi-stable reduction. These groups were studied by Silverberg and Zarhin and more recently by the author in order to give a group theoretic characterization of them depending only on the dimension. We give an overview of the developments since Grothendieck's definition with the added novelty of the case of equal characteristic local fields. 
\end{abstract}

\maketitle

\section{Introduction}

\subsection{The case of elliptic curves}

\subsubsection{} Let $A$ and $B$ be integers such that the curve
\[
E\colon y^2=x^3+Ax+B
\]
is an elliptic curve with discriminant $\Delta_E=-16(4A^3+27B^2)\neq 0$. In order to define the reduction of $E$ at a prime $p$ we can naively consider the curve
\[
E_p\colon y^2=x^3+\overline{A}x+\overline{B}
\]
over $\F_p$. This process depends on the starting equation and without going into technical details, if we assume the starting equation to satisfy a condition of minimality with regards to $p$, then $E_p$ is well defined in this way. Note that since reduction modulo $p$ is a map of rings the discriminant $\Delta_p$ of $E_p$ verifies
\[
\Delta_p= \Delta_E \mod p.
\]
 In particular, if $\Delta_E$ is divisible by $p$ we see that $E_p$ is a singular curve and not an elliptic curve. 

\subsubsection{} There are, in fact, three possibilities for the curve $E_p$.
\begin{itemize}
    \item If the curve $E_p$ is an elliptic curve. We say $E$ has good reduction at $p$.
    \item If the curve $E_p$ has a singular point $x$ and $E_p\setminus\{x\}$ is a torus, we say $E$ has bad multiplicative reduction at $p$.
    \item If the curve $E_p$ has a singular point $x$ and $E_p\setminus \{x\}$ is isomorphic to the additive group $\G_a$, we say $E$ has bad additive reduction at $p$.
\end{itemize}

If $E_p$ falls into one of the first two cases, we also say $E$ has \textit{semi-stable} reduction at $p$. 

\medskip

Our main interest is to understand how these different types of reduction evolve when we take field extensions. Let us showcase this with a concrete example. Consider, over the rationals, the curve $E$ defined by the equation
\[
E\colon y^2=x^3+3p^2x-2p^3,\text{ with discriminant }\Delta_E=-2^7\cdot 3^3\cdot p^6.
\]
One can check that this equation is minimal at $p$ and it is clear that $E$ has bad reduction at $p$. Now let us consider $E$ as a curve over the field $\Q(\sqrt{p})$. By adding $\sqrt{p}$ we have new changes of variables available. Let us make the change
\[
x'=px,~y'=p\sqrt{p}y
\]
which gives
\[
E\colon {y'}^2={x'}^3+3p^2x'-2p^3
\]
and finally 
\[
E\colon y^2= x^3+3x-2.
\]
We can check that this equation is minimal at $p$ over $\Q(\sqrt{p})$ and its discriminant is $-2^7\cdot 3^3$. In other words, $E$ has good reduction at $p$ over the field $\Q(\sqrt{p})$. The obstruction for $E$ to have good reduction at $p$ has been lifted by taking a field extension. It is these types of possible obstructions that we characterize by group theoretic means in the general case of abelian varieties. 

\subsection{Abelian varieties and the obstruction to semi-stable reduction}

\subsubsection{} In our context, an abelian variety is a smooth complete group variety over a field. The abelian varieties of dimension $1$ are precisely the elliptic curves. The previous section generalizes to higher dimension as follows. For an abelian variety $A$ over a number field $K$ with residue field $k$ at a prime $\mathfrak{p}$, the reduction $A_k$ of $A$ at $\mathfrak{p}$ is defined by using the appropriate fiber of its Néron model. The Néron model of an abelian variety over a number field is a geometrical space fulfilling a universal property, called the Néron mapping property. To be precise, it is a smooth group scheme over the spectrum $\Spec \OO_K$ of the ring of integers of $K$. For our purpose, the important point is that, by a result of Chevalley, $A_k$ fits into a short exact sequence 
\[
\begin{tikzcd}
    0 \arrow[r] & T_0\times U_0 \arrow[r] &A_k \arrow[r] & B_0 \arrow[r] & 0
\end{tikzcd}
\]
where $T_0$ is a torus, $U_0$ is a unipotent algebraic group and $B_0$ is an abelian variety. We call the dimension of $T_0$ the toric rank of the reduction of $A$, the dimension of $U_0$ its unipotent rank and the dimension of $B_0$ its abelian rank. The toric and abelian rank are non decreasing when taking field extensions while the unipotent rank is non increasing. 

\medskip

\begin{defin}
    The abelian variety $A$ has semi-stable reduction at $\mathfrak{p}$ if the unipotent rank of the reduction of $A$ at $\mathfrak{p}$ is $0$. We say $A$ has semi-stable reduction if it has semi-stable reduction at all primes.
\end{defin}

One important fact about the semi-stable situation is that it is invariant by field extensions. A fundamental result of Grothendieck asserts that we can always get to this situation. 

\begin{theo}[\cite{sga} exposé IX]
    There is a finite extension $L/K$ such that the base change $A_L$ of $A$ to $L$ has semi-stable reduction. 
\end{theo}

From the theorem we get that the notion of semi-stable toric rank and semi-stable abelian rank of an abelian variety is well-defined. We will use it. 

\medskip

In proving his theorem, Grothendieck shows that the local obstruction at a given prime $\mathfrak{p}$ above $p$ is governed by a finite group which we denote by $\Phi_{A,v}$ where $v$ is the valuation associated to the prime $\mathfrak{p}$. These groups are called the \textit{finite monodromy groups} of the abelian variety $A$. 

\subsubsection{} The goal of the present text is to give an overview of the work from Grothendieck via Silverberg and Zarhin to the author with the aim of characterizing these groups. The central question around which we will develop the different results is as follows. 

\begin{quest} 
    What can be said about the finite groups $G$ that are finite monodromy groups for some abelian variety of dimension $g$ ?
\end{quest}

We will answer this question first in terms of the \emph{structure} of $G$ and its \emph{cardinality} and then by going to \emph{characterization} of such groups. In some cases, $g=1$ or $g=2$, the complete list of such groups is provided. 

\medskip

 The notion of $(p,t,a)$-inertial groups was introduced by Silverberg and Zarhin, see Definition~\ref{def:pta}, for this purpose. The main result can be stated in the following way. 

\begin{theo} \label{theo:main}
     Let $G$ be a finite group. Then $G$ is $(p,t,a)$-inertial if and only if $G$ is the finite monodromy group of an abelian variety $A$ of dimension $t+a$ over a $p$-adic field $K$ (resp. a local field of equal characteristic $K$ with perfect residue field) at a place $v$ of residue characteristic $p$ and $A$ is such that, for any extension $L/K$ for which $A_L$ has semi-stable reduction the reduction of $A_L$ has toric rank $t$ and abelian rank $a$. 
\end{theo}

The case of $p$-adic fields follows from \cite{Ph25}, Corollary~4.9. The case of local fields of equal characteristic will be dealt with in this text. 

\subsubsection{} In Section~\ref{sec:grothtosz} we start by going over the first properties of the groups $\Phi_{A,v}$ that follow from their definition by Grothendieck. Afterwards we go into more precise results by Silverberg and Zarhin in \cite{SZ95, SZ98, SZ05} which lead to the notion of $(p,t,a)$-inertial groups as well as the realization of the list of $(p,t,a)$-inertial groups with $t+a=2$ as the local obstruction for abelian surfaces over equal characteristic local fields.

\medskip

Section~\ref{sec:auth} gives an account of the different characterizations of the finite monodromy groups over number fields by the author in \cite{Phi24, Ph25}. A first characterization, that is geometric in nature, using semi-abelian varieties over finite fields, is obtained by $p$-adic Hodge theory. The main result then follows by a careful study of the rational group algebras of $(p,t,a)$-inertial groups and their representations.

\section{The groups that represent the local obstruction to semi-stable reduction} \label{sec:grothtosz}

We start by recalling Grothendieck's definition of the finite monodromy groups of abelian varieties and move towards their properties studied by Silverberg and Zarhin up to the notion of $(p,t,a)$-inertial groups. We then give parts of the complete list of finite monodromy groups for abelian surfaces as shown by Silverberg and Zarhin and Chrétien-Matignon in the case of number fields.

\subsection{Definition and basic properties}

\subsubsection{} Let $A$ be an abelian variety over a number field $K$ and $v$ a non archimedean place of $K$ of residue characteristic $p$. For $\ell\neq p$ a prime number, the $\ell$-adic Tate module of $A$ carries important arithmetic information. It comes equipped with a representation of the absolute Galois group $G_K$ of $K$. It can be succinctly defined as follows. For $n$ a positive integer, the group of $n$-torsion points $A[n]$ of $A(\overline{K})$ is naturally equipped with maps $A[\ell^{n-1}]\to A[\ell^n]\to A[\ell^{n+1}]$ and with an action of $G_K$ compatible with those maps. We can thus form their inverse limit $\varprojlim_n A[\ell ^n]$ which is by definition the $\ell$-adic Tate module $\operatorname{T}_{\ell} A$ of $A$. By construction, with the identification $\operatorname{T}_{\ell}A\simeq \Z_{\ell}^{2g}$, we get an $\ell$-adic representation
\[
\rho_{A,\ell} \colon G_K\longrightarrow \GL_{2g}(\Q_{\ell})\simeq \Aut( \operatorname{T_{\ell} A)\otimes_{\Z_{\ell}} \Q_{\ell}}.
\]
In \cite{sga} exposé IX, Grothendieck shows that there is a biggest open subgroup $I_{v,A}$ of the absolute inertia group $I_v\subset G_K$ at $v$ such that, for any prime $\ell\neq p$, the action of any element of $I_{v,A}$ on the $\ell$-adic Tate module $T_{\ell} A$ of $A$ is unipotent.

\begin{prop} \label{prop:corpskva}
    There is a smallest extension $K_{v,A}$ of the maximal unramified extension $K_{v}^{\mathrm{un}}$ of $K_v$ such that $A_{K_{v,A}}$ has semi-stable reduction. Furthermore, this extension is finite Galois with group $\Phi_{A,v}=I_v/I_{v,A}$. 
\end{prop}

The groups $\Phi_{A,v}$ obtained in this way are our main object of study.

\begin{defin}The group $\Phi_{A,v}=I_v/I_{v,A}$ is called the \textit{finite monodromy group} of $A$ at $v$. 
\end{defin}

The term finite monodromy group is first used by Raynaud in \cite{Ray99}. An equivalent definition for these groups was given by Silverberg and Zarhin in \cite{SZ95}. 

\begin{theo}[\cite{SZ95}, Theorem~5.2]
    Let $G_v$ be the Zariski closure of the image of $I_v$ under the $\ell$-adic representation
    \[
    \rho_{A,\ell} \colon G_K\to \operatorname{GL}_{2g}(\Q_{\ell}).
    \]
    The group $\Phi_{A,v}$ is isomorphic to the group of components $G_v/G_v^{\circ}$ of $G_v$.
\end{theo}

Let us just note that $G_v/G_v^{\circ}$ is a finite étale algebraic group, but since by definition the $\Q_{\ell}$-rational points of $G_v$ are dense, its component group $G_v/G_v^{\circ}$ can be identified with the finite group $G_v/G_v^{\circ}(\Q_{\ell})$.

\begin{rem}
    In the case of elliptic curves these groups were first introduced by Serre in \cite{Se72} section~5.6 with a list of the possible groups given by :
    \begin{itemize}
        \item For $p=2$ the possible groups are $\Z/2\Z$, $\Z/3\Z$, $\Z/4\Z$, $\Z/6\Z$, $Q_8$ and $\operatorname{SL}_2(\F_3)$. 
        \item  For $p=3$ the possible groups are $\Z/2\Z$, $\Z/3\Z$, $\Z/4\Z$, $\Z/6\Z$ and $\Z/3\Z \rtimes \Z/4\Z$. 
        \item For $p\geq 5$ the possible groups are $\Z/2\Z$, $\Z/3\Z$, $\Z/4\Z$ and $\Z/6\Z$. 
        
    \end{itemize}
\end{rem}

\subsubsection{} A few basic properties of these groups follow directly from their definition. In order to present them we first need to introduce ramification groups.

\begin{defin}
    A finite group $G$ is said to be a ramification group at $p$ if it is the inertia group of a Galois extension of local fields of residue characteristic $p$.
\end{defin}

When the context is clear we will only say that $G$ is a ramification group. By definition, the finite monodromy group $\Phi_{A,v}$ is a ramification group at $p$. 

\begin{prop}
    The group $\Phi_{A,v}$ has the following properties.
    \begin{itemize}
        \item[(i)] The abelian variety $A$ has semi-stable reduction at $v$ if and only if $\Phi_{A,v}=\{1\}$
        \item[(ii)] There is an isomorphism
        \[
        \Phi_{A,v}\simeq \Gamma_p\rtimes \Z/n\Z
        \]
        where $\Gamma_p$ is a $p$-group and $n$ is an integer prime to $p$.
        \item[(iii)] The group $\Phi_{A,v}$ is invariant under isogeny, that is, if $B$ is an abelian variety isogenous to $A$ then $\Phi_{B,v}\simeq \Phi_{A,v}$.
    \end{itemize}
\end{prop}

 The first point follows from Grothendieck's definition. The second is a standard fact about ramification groups and gives a strong constraint on the structure of these groups. The third follows directly from the fact that $\Phi_{A,v}$ is determined by the $\ell$-adic representation $\rho_{A,\ell}$, which is itself an isogeny invariant. 

 \medskip

The first point already gives meaning to the fact that $\Phi_{A,v}$ represents the local obstruction to semi-stable reduction at $v$ for $A$. The following result show that the connection goes deeper. 

\begin{prop}
    There is an extension $L/K$ with $[L:K]=\Card \Phi_{A,v}$ such that $A_L$ has semi-stable reduction at $v$. Moreover, if $L/K$ is an extension and $w\mid v$ is a place of $L$ such that $A_L$ has semi-stable reduction at $w$ then the ramification degree $e(w/v)$ is divisible by $\Card \Phi_{A,v}$.
\end{prop}

One can actually show that, in a sense, the finite monodromy groups also represent the global obstruction to semi-stable reduction.

\begin{prop}[\cite{Phi22}, Theorem~2.4]
    Let $\Sigma_K$ be the set of non-archimedean places of $K$. There is an extension $L/K$ with 
    \[
    [L:K]=\underset{v\in \Sigma_K}{\operatorname{lcm}} \Card \Phi_{A,v}
    \]
    such that $A_L$ has semi-stable reduction. Moreover, if $L/K$ is an extension such that $A_L$ has semi-stable reduction then 
    \[
    \underset{v\in \Sigma_K}{\operatorname{lcm}} \Card \Phi_{A,v} \mid [L:K].
    \]
\end{prop}

\subsection{Restrictions on finite monodromy groups and $(p,t,a)$-inertial groups}

\subsubsection{} Let us first introduce some context. We now consider an abelian variety $A$ over the maximal unramified extension $K$ of a local field and we denote its valuation by $v$. Let $k$ be the residue field of $K$. In order to study the finite monodromy group of an abelian variety over a number field at a given place we can always base change to this situation by invariance of the Néron model under unramified extensions. Let $A'$ be the base change of $A$ to the field $K_{v,A}$ given by Proposition~\ref{prop:corpskva}. By definition $A'$ has semi-stable reduction, that is $A'_k$ is a semi-abelian variety.

\medskip

The following basic idea, already present in \cite{ST68} relates the group $\Phi_{A,v}$ to the reduction $A'_k$ of $A'$. The variety $A'$ carries a descent datum for the extension $K_{v,A}/K$ which corresponds to an action of the Galois group $\Phi_{A,v}$ on $A'$ seen as a $K$-scheme. By the Néron mapping property the action of $\Phi_{A,v}$ on $A'$ extends to its Néron model $\mathcal{A}'$. Let $\sigma\in \Phi_{A,v}$. The pullback to the special fiber of the square
\[
\begin{tikzcd}
\mathcal{A}' \arrow[r] \arrow[d] & \mathcal{A}'  \arrow[d] \\
\Spec \OO_{K_{v,A}} \arrow[r] & \Spec \OO_{K_{v,A}}
\end{tikzcd}
\]
gives a triangle
\[
\begin{tikzcd}
    A'_k \arrow[rr] \arrow[rd] & & A'_k \arrow[ld] \\
    & \Spec k & 
\end{tikzcd}
\]
and thus a map $\Phi_{A,v}\to \Aut (A'_k)$. The main point is that this map is injective. 

\subsubsection{} This idea was then refined in \cite{SZ98}. There is a natural map $A_k \to A'_k$ coming from the base change of the Néron models of $A$ and $A'$. Recall that $A_k$ sits in an exact sequence
\[
\begin{tikzcd}
    0 \arrow[r] & T_0\times U_0 \arrow[r] & A_k \arrow[r] & B_0 \arrow[r] & 0. 
\end{tikzcd}
\]
Let us denote by $t$ the toric rank of $A_k$ and by $a$ its abelian rank. Since $A'$ has semi-stable reduction, it sits in a similar exact sequence but with no unipotent part
\[
\begin{tikzcd}
    0 \arrow[r] & T'_0 \arrow[r] & A_k \arrow[r] & B'_0 \arrow[r] & 0. 
\end{tikzcd}
\]

We denote by $t'$ its toric rank and $a'$ its abelian rank. From the base change map of Néron models we get induced maps $T_0\to T'_0$ and $B_0\to B'_0$. Let us denote by $T$ and $B$ the corresponding quotients $T'_0/T_0$ and $B'_0/B_0$. They are respectively of dimension $t'-t$ and $a'-a$.

\medskip

The main observation is now that the image of the map $A_k\to A'_k$ is fixed by the action of $\Phi_{A,v}$ on $A'_k$ and thus $\Phi_{A,v}$ acts on $T$ and $B$. By further considering that $A'$ can be equipped by a polarization, which is coming from $A$ and thus invariant by the Galois action of $\Phi_{A,v}$, one gets the following result.

\begin{theo}[\cite{SZ98}, Theorem~5.3] \label{theo:szladic}
    Let $\lambda$ be a polarization on $B$ invariant under the action of $\Phi_{A,v}$ and let $\ell$ be a prime that does not divide $\deg \lambda$. Then we have the following injections
    \[
    \Phi_{A,v} \hookrightarrow \Aut T \times \Aut (B,\lambda)
    \]
    and 
    \[
    \Phi_{A,v} \hookrightarrow \GL_{t'-t}(\Z) \times \Sp_{2(a'-a)} (\Z_{\ell}).
    \]
\end{theo}

Such a polarization $\lambda$ always exists on $B$ since $\Phi_{A,v}$ is a finite group. As for the proof, one obtains the second map by passing to the $\ell$-adic Tate module of $B$ from the first and then one shows that it is an injection using the faithful action of $\Phi_{A,v}$ on $\operatorname{T}_{\ell} A'_k\simeq (\operatorname{T}_{\ell} A)^{I_{v,A}}$. 

\subsubsection{} From Theorem~\ref{theo:szladic} one can get numerical bounds on the order of $\Phi_{A,v}$. To that end, let us introduce, for $n$ a positive integer, 
\[
M(n)=\operatorname{lcm} \{\Card G\mid G\subset \GL_n(\Q),~G \text{ finite} \}.
\] 
This function is called the Minkowski bound as it was computed by Minkowski in 1887 to be
\[
M(n)= \prod\limits_p p^{r(n,p)} \text{ where } r(n,p)=\sum\limits_{i\geq 0} \Big\lfloor \frac{n}{p^i(p-1)} \Big\rfloor,\text{ for any integer } n. 
\]

\begin{cor}[\cite{SZ98}, Corollary~6.3]
    We have the following sequence of divisibilities 
    \[
    \Card \Phi_{A,v} \mid M(t'-t)\cdot M(2(a'-a)) \mid M(2g).
    \]
\end{cor}

\subsubsection{} \label{subsub:pta} Finally, we recall the notion of $(p,t,a)$-inertial groups from \cite{SZ05}. 

\begin{defin}[\cite{SZ05} Definition 1.1] \label{def:pta} Let $p$ be a prime number or $p=0$ and $t$, $a$ positive integers. A finite group $G$ is said to be $(p,t,a)$-inertial if it satisfies the two following conditions :
\begin{itemize}
    \item[$(i)$] If $p=0$ then $G$ is cyclic otherwise $G$ is a semi-product $\Gamma_p\rtimes \Z/n\Z$ with $\Gamma_p$ a $p$-group and $n$ an integer prime to $p$.
    \item[$(ii)$] For all primes $\ell\neq p$ there is an injection
    \[
    \iota_{\ell}\colon G\hookrightarrow \GL_t(\Z) \times \Sp_{2a} (\Q_{\ell})
    \]
    
\end{itemize}
    such that the projection map onto the first factor is independent of $\ell$ and the characteristic polynomial of the projection of any element onto the second factor has integer coefficients independent of $\ell$. 
\end{defin}

Let us remark that, from the last part of the definition, the character $\chi_{\ell}$ of the representation of $G$ given by the projection to the second factor is independent of $\ell$. 

\medskip

The purpose of this notion is to characterize finite monodromy groups. A precise statement is given by Question~1.13 of \cite{SZ05}. A stronger, but essentially equivalent question is as follows. \emph{Let $G$ be a $(p,t,a)$-inertial group. Is $G$ the finite monodromy group of an abelian variety $A$ over a $p$-adic, or local field of equal characteristic $p$, for which the toric rank of the semi-stable reduction of $A$ is $t$ its abelian rank is $a$ ?}

\medskip

One immediate property of these groups is that if $G$ is $(p,t,a)$-inertial then for any $s\geq t$ and any $b\geq a$ the group $G$ is also $(p,s,b)$-inertial. It is clear that an analogous property for finite monodromy groups holds.

\medskip

From the previous paragraph we get that the group $\Phi_{A,v}$ is $(p,t,a)$-inertial with $p$ the residue characteristic of $v$, $t$ the toric rank of $A_k$ and $a$ its abelian rank. More specifically, $\Phi_{A,v}$ is $(p,t'-t,a'-a)$-inertial. One can show even more. For a $(p,t,a)$-inertial group $G$ and a given family of maps $(\iota_{\ell}\colon G\to \GL_t(\Z)\times \Sp_{2a}(\Q_{\ell}))_{\ell\neq p}$ we can consider the rank of the spaces $(\Z^t)^G$ and $(\Q_{\ell}^{2a})^G$ of $G$-invariants. The rank of the space $(\Q_{\ell}^{2a})^G$ is independent of $\ell$ by the assumptions on the family $(\iota_{\ell})_{\ell\neq p}$. 

\begin{defin}[\cite{SZ05} Definition~1.10]
    A $(p,t,a)$-inertial group $G$ is said to be strongly $(p,t,a)$-inertial if there is a family of maps $(\iota_{\ell}\colon G\to \GL_t(\Z)\times \Sp_{2a}(\Q_{\ell}))_{\ell\neq p}$ satisfying condition $(ii)$ of Definition~\ref{def:pta} and such that 
    \[
     \operatorname{rank} (\Z^t)^G= \operatorname{rank} (\Q_{\ell}^{2a})^G=0.
    \]
\end{defin}

By the ideas presented in the previous paragraphs one can show that the finite monodromy group $\Phi_{A,v}$ is strongly $(p,t'-t,a'-a)$-inertial. 

\medskip

We finish by going back to Question~1.13 of \cite{SZ05}. If $G$ is $(p,t,a)$-inertial but not strongly $(p,t,a)$-inertial, then we can quotient out the multiples of the trivial representation present to turn $G$ into a strongly $(p,t',a')$-inertial group with $t'\leq t$ and $a'\leq a$. Then, if $G$ is realized as a finite monodromy group of an abelian variety $A$ over $K$ with semi-stable toric rank $t'$ and semi-stable abelian rank $a'$, we can give a positive answer to the stronger question by considering $A\times B$, where $B$ is a semi-stable abelian variety over $K$ with toric rank $t-t'$ and abelian rank $a-a'$. That is, to answer positively Question~1.13 of \cite{SZ05}, it is equivalent to answer positively the question asked in paragraph~\ref{subsub:pta}. 

\subsection{The case of abelian surfaces} \label{sub:galoistwist}

\subsubsection{} Using group theoretic arguments, Silverberg and Zarhin in \cite{SZ05} compute the complete list of $(p,t,a)$-inertial groups when $t+a= 2$. We will only give the generic portions of the list and refer to Definition~1.2 of \cite{SZ05} for the rest. Let us denote by $\Sigma_p(t,a)$ the set of $(p,t,a)$-inertial groups. Then we have

\begin{itemize}
    \item $\Sigma_p(2,0)=\Sigma_p(1,1)= \{\Z/2\Z, \Z/3\Z, \Z/4\Z, \Z/6\Z\}$ for $p\geq 5$;
    \item $\Sigma_p(0,2)=\Sigma_p(2,0) \cup \{ \Z/5\Z, \Z/8\Z, \Z/10\Z, \Z/12\Z\}$ for $p\geq 7$.
\end{itemize}

To realize these groups as finite monodromy groups over local field of equal characteristic the main tool is the \emph{Galois twisting construction} which we briefly recall in our context. 

\begin{defin} For an abelian variety $A$ over a field $K$ and a given finite Galois extension $L/K$ we say that $B$ is an $L/K$-twist of $A$ if there is an isomorphism $B_L\simeq A_L$. 
\end{defin}

The main classification result for twisted abelian varieties is that there is a bijection between such twists $B$ of a given abelian variety $A$ over $K$ and the pointed set in Galois cohomology $H^1(\Gal (L/K), \Aut (A_L))$. In particular, when $\Aut (A_L)$ is a trivial $\Gal (L/K)$-module this set is given by the homomorphisms $\Gal(L/K) \to \Aut (A_L)$. In the case where $A$ has semi-stable reduction the finite monodromy group of the twisted variety $B$ at a given place $v$ can be expressed as the inertia group of $L/K$ at $v$. 

\begin{theo}[\cite{SZ05}, Theorem~4.3] \label{theo:galoistwist}
    Let $A$ be an abelian variety over $K$ with semi-stable reduction at $v$. Let $L/K$ be a Galois extension with group $G$ such that $\Aut (A_L)$ is a trivial $G$-module. Then the twisted variety $B$ corresponding to an injective morphism
    \[
     G\hookrightarrow \Aut A_L
    \]
    verifies that $\Phi_{v,B}$ is the inertia subgroup of $G$. 
\end{theo}

We can thus use this technique to produce abelian varieties with a given finite monodromy group. It requires two main ingredients. First a result in inverse Galois theory. 

\begin{theo}[\cite{SZ05}, Lemma~5.1] If $k$ is an algebraically closed field of characteristic $p$ then every $(p,t,a)$-inertial group can be realized as a Galois group of a totally ramified extension of $k((t))$. 
\end{theo}

Secondly, we need a starting object to do the twisting. That is, in our case, an abelian variety of dimension $t+a$ which has semi-stable reduction and has a subgroup of its automorphism group isomorphic to the given $(p,t,a)$-inertial group we want to realize. In the case of equal characteristic, this last constraint is not very severe by considering a base change of an appropriate abelian surface over a finite field. Such abelian surfaces are given by ad hoc constructions in the case of $G$ being a $(p,t,a)$-inertial group with $t+a=2$. The end result, putting together the different parts, can be stated as follows.

\begin{theo}[\cite{SZ05}]
    Let $G$ be a $(p,t,a)$-inertial group with $t+a=2$. Then there is an abelian variety $A$ over a local field of equal characteristic with valuation $v$ such that
    \[
    \Phi_{A,v}=G.
    \]
\end{theo}

Let us also note that the list was also realized over mixed characteristic local fields. The last missing group was realized by the ad hoc construction of a family of hyperelliptic curves by Chrétien and Matignon in \cite{CM13}. 

\medskip

In \cite{Phi22b}, the author used a variation on this technique to provide abelian varieties over number fields with large finite monodromy groups. To be precise, the finite monodromy groups produced are given by $p$-Sylow subgroups of some matrix groups used by Minkowski and they attain the value of his bound. 

\section{The finite monodromy groups characterized as $(p,t,a)$-inertial groups} \label{sec:auth}

In this section we present recent work by the author that gives a positive answer to the question by Silverberg and Zarhin recalled in paragraph~\ref{subsub:pta}. We start by presenting a geometric characterization of finite monodromy groups for $p$-adic fields before going back to $(p,t,a)$-inertial groups and dealing with the case of local fields of equal characteristic.     

\subsection{A first characterization of finite monodromy}

\subsubsection{}  We work over a $p$-adic field $K$ with residue field $k$. This characterization is given in terms of a descent datum from a variety with semi-stable reduction. The advantage with this characterization is that we have much more control and tools available to work with varieties which have semi-stable reduction. In particular, we can produce such semi-stable abelian varieties by some deformation process. 

\begin{prop}[cf. Theorem~3.9 of \cite{Phi24}] \label{prop:chardescent}
    Let $A$ be an abelian variety over $K^{\mathrm{un}}$ and $G$ a finite group. Then $G$ is the finite monodromy group of $A$ if and only there is a Galois extension $L/K^{\mathrm{un}}$ with group $G$ such that $A_L$ has semi-stable reduction and the canonical descent datum of $A_L$ induces an injection $G\hookrightarrow \Aut (A_L)_{\overline{k}}$.
\end{prop}

From the previous section we have already seen how a descent datum relative to a totally ramified extension induces a group homomorphism into the automorphism group of the reduction.

\subsubsection{} Let us first introduce polarizations of semi-abelian varieties over finite fields. 

\begin{defin}\label{def:vsapol}
A morphism $\lambda_0\colon A_0\to A_0^t$ of semi-abelian varieties over a finite field $k$ 
\[
\begin{tikzcd}
    0 \arrow[r] & T_0 \arrow[r] \arrow[d, "\lambda_{T_0}"] & A_0 \arrow[r, "p"] \arrow[d, "\lambda_0"] & B_0 \arrow[r]  \arrow[d, "\lambda_{B_0}"]& 0\\
    0 \arrow[r] & T_0^t \arrow[r] & A_0^t \arrow[r, "p"] & B_0^{\vee} \arrow[r] & 0
\end{tikzcd}
\]
is a polarization if the induced morphism $\lambda_{T_0}$ is an isogeny and the induced morphism $\lambda_{B_0}$ is a polarization. 
\end{defin}

The first characterization of finite monodromy over $p$-adic fields is now as follows.

\begin{theo}[\cite{Phi24}, Theorem~1.1] \label{theo:hodge}
    Let $G$ be a finite group which is a ramification group. Then there is an abelian variety $A$ of dimension $g$ over a $p$-adic field $K$ such that $G$ is the finite monodromy group of $A$ if and only if there is a polarized semi-abelian variety $(A_0,\lambda_0)$ of dimension $g$ over a finite field of characteristic $p$ with an injective map
    \[
    G\hookrightarrow \Aut (A_0,\lambda_0).
    \]
\end{theo}

The proof is done by constructing $A$ through deformation, degeneration and descent. By a result of Serre and Tate, later generalized by Bertapelle and Mazzari in \cite{BM19}, the deformation theory of semi-abelian varieties is given by the one of their $p$-divisible groups. We are thus lead to deform $A_0[p^{\infty}]$, the $p$-divisible group associated to $A_0$, with its given action of $G$. The description of such objects in terms of semi-linear categories is done by $p$-adic Hodge theory. The degeneration step is done through the theory of Faltings and Chai, see \cite{FC90}. In both steps our goal is to keep track of the action of $G$ and turn it into a descent datum in order to apply the characterization of Proposition~\ref{prop:chardescent}. 

\begin{rem}
    One can obtain the analogous result for local fields of equal characteristic by using the Galois twisting construction of Section~\ref{sub:galoistwist} instead of $p$-adic Hodge theory. 
\end{rem}

\subsubsection{} In order to discuss the proof we give a quick account of the parts of $p$-adic Hodge theory that are needed. One goal of $p$-adic Hodge theory is to understand the categories $\operatorname{Rep}_{G_K} (\Z_p)$ and $\operatorname{Rep}_{G_K} (\Q_p)$ of continuous representations of $G_K$ over $\Q_p$ or $\Z_p$. We will only be interested in finite dimensional representations as in our case they will come from the action of $G_K$ on abelian varieties and their Tate modules or $p$-divisible groups. Furthermore, we will consider the subcategories $\operatorname{Rep}^{\mathrm{st}}_{G_K} (\Q_p)$ and $\operatorname{Rep}^{\mathrm{st}}_{G_K} (\Z_p)$ of \emph{semi-stable} such representations. The precise definition of those require a lot of work and the introduction of a certain period ring $B_{\mathrm{st}}$, but we should note that $G_K$-representations coming from the $p$-divisible groups of semi-stable abelian varieties are semi-stable. 

\medskip

Let $K_0\subset K$ be the maximal unramified subfield of $K$. The extension $K_0/\Q_p$ is Galois with cyclic group generated by the Frobenius element $\sigma$. Fontaine shows that the category $\operatorname{Rep}^{\mathrm{st}}_{G_K} (\Q_p)$ is equivalent to the category of admissible filtered $\varphi$-modules $\operatorname{MF}_K^{\varphi,N}$. The objects of $\operatorname{MF}_K^{\varphi}$ are finite dimensional $K_0$-vector spaces equipped with a $\sigma$-semi-linear and invertible map $\varphi$, a nilpotent linear map $N$ such that $N\varphi=p\varphi N$ and a filtration $F$ for which the so-called  admissibility condition holds -- which is quite technical but essentially corresponds to a compatibility between sub-$\varphi$-modules and the filtration $F$. We will in fact only deal with such modules where $N=0$, that is the subcategory $\operatorname{MF}^{\varphi}_K$ of filtered $\varphi$-modules. Starting with $D\in \operatorname{MF}^{\varphi}_K$ and forgetting the filtration, one recovers a $\varphi$-module, an object of $\operatorname{MF}^{\varphi}$ which is itself the isogeny category of $p$-divisible groups over the residue field $k$.

\medskip

In order to deform our starting variety $A_0$ we have to deal with the more complicated category of $p$-divisible groups over $\OO_K$. This one is shown to be equivalent to the category of Breuil-Kisin modules $\operatorname{BT}_{\mathfrak{S}}^{\varphi}$ which are, roughly speaking, free modules over the ring $\mathfrak{S}=\OO_{K_0}[[u]]$ of formal power series over the integers of $K_0$ equipped with a \emph{Frobenius}, that is a $\sigma$-semi-linear map with added constraints. This category is connected to $\varphi$-modules through the special fiber functor, which corresponds algebraically to taking $u=0$ and inverting $p$. More importantly it is related to the category $\operatorname{MF}_K^{\varphi}$ of filtered $\varphi$-modules. Indeed, filtered $\varphi$-modules are the isogeny category for Breuil-Kisin modules as given by Proposition~2.2.2 of \cite{Kis06}. 

\subsubsection{} The different steps of the proof of Theorem~\ref{theo:hodge} are summarized in the diagram given by Figure~\ref{fig:diagramme preuve}. As stated, the goal is to deform $A_0$ in a way that keeps track of the action of $G$. We are thus lead to deform the $p$-divisible group of $A_0$ and, since our problem is invariant by isogeny, we travel through the relevant isogeny categories.

  \usetikzlibrary{decorations.pathmorphing}
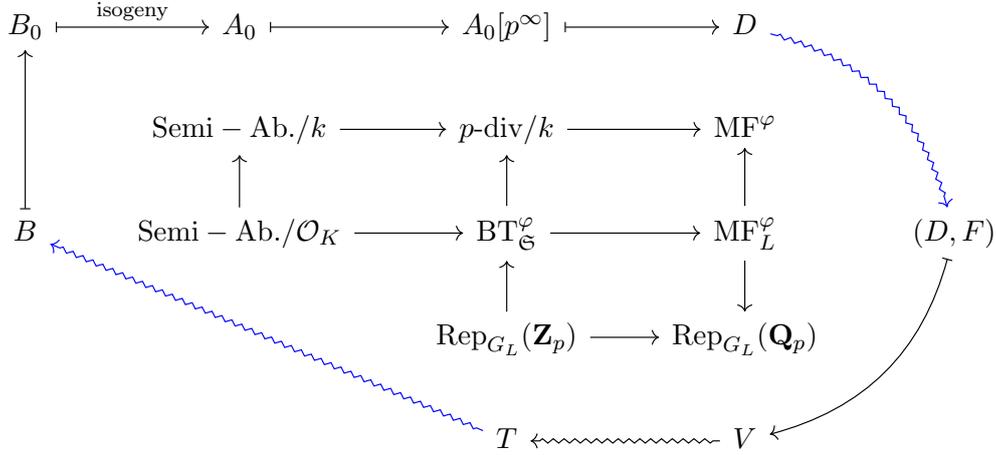
\begin{figure}[H]
    \centering
    \begin{tikzcd}
B_0 \arrow[r, mapsto, "\text{isogeny}"] & A_0 \arrow[r, mapsto] & A_0[p^{\infty}] \arrow[r, mapsto] & D \arrow[ddr, bend left, squiggly, color=blue]\\
& \mathrm{Semi-Ab.}/k \arrow[r] & p\text{-div}/k   \arrow[r] & \operatorname{MF}^{\varphi} & \\
B \arrow[uu, mapsto] & \mathrm{Semi-Ab.}/\OO_K \arrow[r] \arrow[u]& \operatorname{BT}_{\mfS}^{\varphi} \arrow[r] \arrow[u] & \operatorname{MF}_L^{\varphi} \arrow[u] \arrow[d] &  (D,F) \arrow[ddl, bend left, mapsto]\\
& & \operatorname{Rep}_{G_L} (\Z_p) \arrow[r] \arrow[u] & \operatorname{Rep}_{G_L} (\Q_p) & \\
& & T \arrow[uull, squiggly, color=blue] & V \arrow[l, squiggly] & 
\end{tikzcd}
    \caption{Navigating through the categories in play}
    \label{fig:diagramme preuve}
\end{figure}

The blue arrows represent the most difficult steps. The first one, on the right, is about finding a suitable filtration on the Fontaine $\varphi$-module associated to the $p$-divisible group $A_0[p^{\infty}]$. The space of such filtrations is, in general, not a Zariski open subspace of a Grassmannian but only a compact Berkovich open subset. The added difficulty here is that we are looking for a filtration $F$ such that the resulting object $(D,F)$ is equipped with a descent datum that, through the forgetful functor, recovers the action of the group $G$ on $D$. 

\medskip

The second difficult step is to recover an abelian variety as the generic fiber of the  semi-abelian scheme resulting from the deformation process. This is done using Faltings' and Chai's theory of degeneration of abelian varieties. We only state the result here. 

\begin{theo}[\cite{Phi24} Proposition~3.8] \label{theo:faltingschai}
    Let $C$ be a semi-abelian scheme over $\OO_L$ with split toric part and a polarization $\lambda$ which induces an isomorphism on the toric part. Let $(f_{\sigma})_{\sigma\in G}$ be a descent datum on $C$ compatible with $\lambda$ and such that the induced map $\iota\colon G\to \Aut (C_k, \lambda_k)$ is injective. Then there is an abelian variety $A$ over $L$ whose reduction $A_k$ is isomorphic to $C_k$ and which is equipped with a descent datum $(h_{\sigma})_{\sigma\in G}$ that induces the map $\iota$ and which has a polarization of degree $\deg \lambda$.
\end{theo}

Obtaining Theorem~\ref{theo:hodge} from this point is simply an application of Proposition~\ref{prop:chardescent}.

\subsection{Geometric representations of the rational group algebras of ramification groups}

\subsubsection{} The rational group algebra of a ramification group $G$ has a specific structure. A theorem of Serre tells us that $\Q[G]$ is quasi-split outside of $p$, that is, for any prime $\ell\neq p$, $\Q[G]\otimes_{\Q} \Q_{\ell}$ is a product of matrix algebras over fields. We give a refinement of this result which deals with the center of the simple factors of $\Q[G]$.

\begin{theo}[] \label{theo:centresalgebresintro}
   Let $G$ be a ramification group at $p$. The rational group algebra $\Q[G]$ has the following properties. 
    \begin{itemize}
        \item[$(i)$] The algebra $\Q[G]$ is quasi-split outside $p$.
        \item[$(ii)$] If $E\subset \Q[G]$ is a simple factor of $\Q[G]$ then the center $Z(E)$ of $E$ is a CM-field or a subfield of $\Q(\mu_{p^{\infty}})$.
    \end{itemize}
    Furthermore, if $\Card G$ is odd then no simple factor $E$ of $\Q[G]$ is such that $E\otimes_{\Q} \R$ is a product of matrix algebras over the quaternions.  
\end{theo}

One should further note that $\Q[G]$, and thus all its simple factors, is equipped with a natural positive involution given by $g\mapsto g^{-1}$. 

\medskip

Such algebras, of finite dimension over $\Q$ and admitting positive involutions, were classified by Albert in four types, see for example \cite{Mu08} p. 186 (202), which we briefly recall here. The classification is given for skew-fields $D$ with a positive involution and extend to the general case.

\begin{itemize}
    \item Type I : The skew-field $D$ is a totally real number field.
    \item Type II : The skew-field $D$ is a quaternion algebra over its center $F$. The number field $F$ is totally real and for every embedding $\sigma\colon F\to \R$ there is an isomorphism
    \[
    \R\otimes_{F} D \simeq  M_2(\R). 
    \]
    \item Type III : The skew-field $D$ is a quaternion algebra over its center $F$. The number field $F$ is totally real and for every embedding $\sigma\colon F\to \R$ there is an isomorphism
    \[
    \R \otimes_{F} D \simeq \HH
    \]
    where $\HH$ is the Hamilton quaternion algebra over $\R$.
    \item Type IV : The center $F$ of $D$ is a totally imaginary quadratic extension of a totally real field $F_0$ with conjugation $\overline{.}$ over $F_0$. The local invariants from the Brauer-Hasse-Noether theorem are such that for any finite place $v$ such that $\overline{v}=v$ then $\operatorname{inv}_v(D)=0$ and otherwise $\operatorname{inv}_v(D)+\operatorname{inv}_{\overline{v}}(D)=0$. 
\end{itemize}

\subsubsection{} In order to relate $G$ to semi-abelian varieties we show that the simple factors of $\Q[G]$ can be embedded in a nice way into endomorphism algebras of abelian varieties by the study of their representations and the use of Honda-Tate theory.  

\medskip

By the definition of $(p,t,a)$-inertial groups we are interested in representations of simple algebras which are \textit{rational}, that is with characteristic polynomials with rational coefficients and which are \textit{polarized}, that is, which admit an invariant non degenerate alternate bilinear form. Given such a representation of dimension $a$ from a simple factor of $\Q[G]$ we get maps $G\to \Sp_{2a}(\Q_{\ell})$ for $\ell\neq p$ which satisfy the condition $(ii)$ of Definition~\ref{def:pta}. For such representations, the classification result over $\C$ is as follows. 

\begin{theo}[\cite{Ph25} Theorem~2.14] \label{theo:repratetpol}
  Let $E$ be a simple finite dimensional $\Q$-algebra. All the rational (resp. and polarized) $\C$-representations of $E$ are multiples of a unique such representation $V_r(E)$ (resp. a unique such representation $V_{rp}(E)$). Furthermore we have that $\dim_{\C} V_r(E)=\deg E$, $V_{rp}(E)\simeq V(E)$ if $E$ is of type III or IV, and $V_{rp}(E)\simeq V_r(E)^2$ if E is of type I or II.
\end{theo}

We can now state the notion of embeddings we are interested in and the associated result on ramification groups. 

\begin{defin}[\cite{Ph25} Definition~2.15]\label{def:goodembed}
    Let $E$ and $E'$ be polarizable simple $\Q$-algebras. A good embedding is an injective morphism $E\hookrightarrow E'$ such that $V_{rp}(E)=V_{rp}(E')$.
\end{defin}

\begin{theo}[cf. \cite{Ph25} Proposition~2.21 and Theorem~2.22]\label{theo:goodembgeo}
    Let $G$ be a ramification group at $p$ and $E$ a simple factor of $\Q[G]$. There is an abelian variety $A$ over $\overline{\F_p}$ and a good embedding
    \[
    E\longrightarrow \End A \otimes \Q.
    \]
\end{theo}

The abelian variety $A$ is constructed by different means whether $E$ is of type III or not. In case of $E$ of type III it can be shown that we can choose $A$ to be a power of a supersingular elliptic curve. In all other cases $A$ is constructed using Honda-Tate theory. This theory classifies simple abelian varieties up to isogeny over finite fields by algebraic integers called Weil $q$-integers. Here one is able to construct a Weil $q$-integer with the relevant properties by algebraic number theory.

\subsubsection{} To conclude we argue differently depending on whether we are dealing with mixed or equal characteristic. For the mixed characteristic case we use Theorem~\ref{theo:hodge} and for the equal characteristic case we use the Galois twisting construction on the semi-abelian variety $A\otimes_k K$ and apply the degeneration step following Faltings and Chai's theory afterwards. 

\begin{theo}\label{theo:realmonod}
    Let $G$ be a ramification group at $p$ and let $(\iota_{\ell}\colon G\to \GL_t(\Z)\times \Sp_{2a}(\Q_{\ell}))_{\ell\neq p}$ be a family of maps satisfying condition $(ii)$ of Definition~\ref{def:pta}. Then there is an abelian variety $A$ of dimension $t+a$ over a $p$-adic field $K$ (resp. over a local field of equal characteristic with perfect residue field) such that $G$ is the finite monodromy group of $A$. Moreover, the action of $G$ on the reduction $A_0$ of $A_{L_A}$ recovers the family of maps $(\iota_{\ell}\colon G\to \GL_t(\Z)\times \Sp_{2a}(\Q_{\ell}))_{\ell\neq p}$ by passing to the $\ell$-adic Tate modules. 
\end{theo}

We give a sketch of proof in the case of local fields of equal characteristic. The case of $p$-adic fields is given by Theorem~4.8 \cite{Ph25}.

\begin{proof}(sketch)
Let us consider a ramification group $G$ with a family of maps $(\iota_{\ell}\colon G\to \GL_t(\Z)\times \Sp_{2a}(\Q_{\ell}))_{\ell\neq p}$ satisfying condition $(ii)$ of Definition~\ref{def:pta}. 

\medskip

Fixing $\ell$ and considering the $\Q$-algebra generated by the image of $G$ by the second projection into $\mathrm{M}_{2a} (\Q_{\ell})$ we recover a map 
\[
\Q[G] \longrightarrow E
\]
where $E$ is a quotient of $\Q[G]$ equipped with an action on a $\Q_{\ell}$-vector space with an invariant non-degenerate bilinear alternating form. Using Theorem~\ref{theo:goodembgeo} we are able to find an abelian variety $A$ over $\overline{\F_p}$ with a good embedding
\[
E\longrightarrow \End A.
\]
We must furthermore have that $\dim A$ divides $a$ from Tate's results. So, up to replacing $A$ with some power of itself and using a diagonal embedding of $E$, we can assume $\dim A=a$. 

\medskip

We can now consider $A_0=\G_m^t\times A$, a semi-abelian variety over some finite field $k$ of residue characteristic $p$ such that
\[
G \hookrightarrow \Aut A_0\simeq \GL_t(\Z) \times \Aut A
\]
and we recover $\iota_{\ell}$ by going to the automorphisms of the $\ell$-adic Tate module of $A$.

\medskip

In order to finish the proof, we consider a local field $K$ of characteristic $p$ with perfect residue field $k$ which admits a Galois extension $L/K$ of group $G$. By twisting $A_0\times_k L$ using Theorem~\ref{theo:galoistwist} we obtain a semi-abelian variety $B$ over $K$ such that $B_L$, or its Néron model, is equipped with a descent datum verifying the conditions of Theorem~\ref{theo:faltingschai}. The resulting abelian variety over $L$ descends to an abelian variety over $K$ with $G$ as finite monodromy group. The remaining parts of the statement follow by construction.
    
\end{proof}

The main corollary is that we can characterize finite monodromy groups over local fields as $(p,t,a)$-inertial groups as given by Theorem~\ref{theo:main}.

\printbibliography

\end{document}